\documentclass[12pt]{article}
\usepackage{latexsym}
\usepackage{amsfonts}
\usepackage{amssymb}
\usepackage{amscd}
\usepackage{array}
\usepackage{amsmath}
\usepackage{epsfig}
\usepackage{a4}

\begin{document}
\title{\textbf{Stringy Hodge numbers for a class of isolated singularities and for threefolds}\\ \date{}}
\author{Jan Schepers and Willem Veys\footnote{Partially supported by project G.0318.06 from the Research Foundation - Flanders (FWO). During the preparation of the manuscript the first author was a Research Assistant of the Research Foundation - Flanders (FWO).}}

\maketitle
\begin{center}
\footnotesize{\textbf{Abstract}}
\end{center}
{\footnotesize Batyrev has defined the stringy $E$-function for complex varieties with at most log terminal singularities. It is a rational function in two variables if the singularities are Gorenstein. Furthermore, if the variety is projective and its stringy $E$-function is a polynomial, Batyrev defined its stringy Hodge numbers essentially as the coefficients of this $E$-function, generalizing the usual notion of Hodge numbers of a nonsingular projective variety. He conjectured that they are nonnegative. We prove this for a class of `mild' isolated singularities (the allowed singularities depend on the dimension). As a corollary we obtain a proof of Batyrev's conjecture for threefolds in full generality. In these cases, we also give an explicit description of the stringy Hodge numbers and we suggest a possible generalized definition of stringy Hodge numbers if the $E$-function is not a polynomial.}
 \\
\section{Introduction}

\noindent \textbf{1.1.} The stringy $E$-function is an interesting singularity invariant, introduced by Batyrev. In \cite{Batyrev} it is associated to complex algebraic varieties with at most log terminal singularities, and it is for instance used to formulate a topological mirror symmetry test for Calabi-Yau varieties with singularities. 

In order to define this invariant and to recall an intriguing conjecture of Batyrev, we quickly review Hodge-Deligne polynomials and some birational geometry.\\

\noindent \textbf{1.2.} We work over the base field of complex numbers; a variety is assumed to be irreducible. We will denote the (singular) cohomology and compactly supported cohomology with coefficients in $\mathbb{C}$ of an algebraic set $X$ by $H^{\bullet}(X)$ and $H^{\bullet}_c(X)$, respectively.

When $X$ is smooth and projective, its cohomology $H^i(X)=H^i_c(X)$ carries a pure Hodge structure of weight $i$. Its Hodge numbers $h^{p,q}(X)$ are given by the dimensions of the $(p,q)$-Hodge component of $H^{p+q}(X)$. They satisfy the following symmetry and Serre duality:
\[ h^{p,q}(X)=h^{q,p}(X)=h^{d-p,d-q}(X),\tag{$1$}  \]
where $d$ is the dimension of $X$.

For an arbitrary $X$ we have that $H^i_c(X)$ carries a mixed Hodge structure, see for example \cite{Deligne1}, \cite{Deligne2}, \cite{Srinivas}. We denote its $(p,q)$-component by $H^{p,q}(H^i_c(X))$ and the dimension of the latter by $h^{p,q}(H^i_c(X))$. Note that here also for other $i$, besides $p+q$, we can have a nontrivial $(p,q)$-component. One defines the Hodge-Deligne polynomial $H(X;u,v)\in \mathbb{Z}[u,v]$ of $X$ as
\[H(X)=H(X;u,v):=\sum_{p,q=0}^d \left[\sum_{i=0}^{2d} (-1)^i h^{p,q}(H^i_c(X))\right] u^pv^q,\]
where again $d=\dim X$. It is a generalized Euler characteristic; so for a closed algebraic subset $Y\subset X$, $H(X)=H(Y)+H(X\setminus Y)$, and for algebraic sets $X$ and $X'$, $H(X\times X')=H(X)\cdot H(X')$. Returning to a smooth and projective $X$, the symmetries and dualities (1) are translated in terms of $H(X)$ as 
\[H(X;u,v)=H(X;v,u) \quad \text{and} \quad H(X;u,v)=(uv)^dH(X;u^{-1},v^{-1}). \tag{$2$}\]
The first equality is still true for arbitrary $X$, but the second in general \textsl{not}.\\

\noindent \textbf{1.3.} A normal variety $X$ is called Gorenstein if its canonical divisor $K_X$ is Cartier, and $\mathbb{Q}$-Gorenstein if $rK_X$ is Cartier for some $r\in\mathbb{Z}_{>0}$. (Note that for instance all hypersurfaces, and more generally all complete intersections are Gorenstein.)

Let $X$ be $\mathbb{Q}$-Gorenstein. Take a log resolution $\varphi:\widetilde{X}\to X$ of $X$, i.e. $\varphi$ is a proper birational morphism from a smooth variety $\widetilde{X}$ such that the exceptional locus of $\varphi$ is a divisor whose irreducible components $D_i$ are smooth and have normal crossings. Then we have $rK_{\widetilde{X}}-\varphi^*(rK_X)=\sum_i b_iD_i$, with $b_i\in \mathbb{Z}$. Usually this is formally rewritten as $K_{\widetilde{X}}-\varphi^*(K_X)=\sum_i a_iD_i$, where $a_i=\frac{b_i}{r}$. The divisor $K_{\widetilde{X}|X}:=K_{\widetilde{X}}-\varphi^*(K_X)$ is called the discrepancy divisor of $\varphi$, and $a_i$ the discrepancy of $D_i$. 

The $\mathbb{Q}$-Gorenstein variety $X$ is called terminal, canonical and log terminal if $a_i>0$, $a_i\geq 0$ and $a_i>-1$ for all $i$, respectively (this is independent of the chosen resolution). Such varieties have `mild' singularities and show up in the Minimal Model Program.\\

\noindent \textbf{1.4.} \textbf{Definition} \cite{Batyrev}\textbf{.} Let $X$ be a log terminal variety. Take a log resolution $\varphi: \widetilde{X}\to X$ with irreducible components $D_i, i\in T,$ of the exceptional divisor. For $I\subset T$, denote $D_I:=\cap_{i\in T} D_i$ and $D_I^{\circ}:=D_I\setminus (\cup_{j\in T\setminus I} D_j)$; so in particular $D_{\emptyset}=\widetilde{X}$. We have the standard stratification $\widetilde{X}=\coprod_{I\subset T} D_I^{\circ}$ into locally closed subsets. The \textsl{stringy $E$-function} of $X$ is
\[E_{st}(X)=E_{st}(X;u,v):=\sum_{I\subset T} H(D_I^{\circ};u,v) \prod_{i\in I} \frac{uv-1}{(uv)^{a_i+1}-1},\] 
where $a_i$ is the discrepancy of $E_i$. Batyrev proved that this expression is independent of the chosen log resolution by means of motivic integration. (Alternatively, one can use the Weak Factorization Theorem of Abramovich, Karu, Matsuki and W{\l}odarczyk from \cite{weakfactorization}.) Concrete computations of the stringy $E$-function in special cases can be found in \cite{Dais} and \cite{DaisRoczen}, but pay attention to the errors in the latter, as explained in \cite{schepers}.\\

\noindent \textbf{Remark.} The following first three statements are clear. The fourth is an easy calculation \cite{Batyrev}.
\begin{description} 
\item[(i)] If $X$ is smooth, then $E_{st}(X)=H(X)$. If $X$ admits a crepant resolution $\varphi:\widetilde{X}\to X$, i.e. such that $K_{\widetilde{X}|X}=0$, then $E_{st}(X)=H(\widetilde{X})$.
\item[(ii)] If $X$ is Gorenstein we have that all $a_i\in \mathbb{Z}_{>0}$. Then $E_{st}(X)$ is a rational function in $u$ and $v$, more precisely it is an element of $\mathbb{Z}[[u,v]]\cap \mathbb{Q}(u,v)$. 
\item[(iii)] Just as Hodge-Deligne polynomials, stringy $E$-functions are symmetric in $u$ and $v$: $E_{st}(X;u,v)=E_{st}(X;v,u)$. 
\item[(iv)] An alternative expression for $E_{st}(X)$ is 
\[E_{st}(X;u,v)= \sum_{I\subset T} H(D_I;u,v)\prod_{i\in I} \frac{uv-(uv)^{a_i+1}}{(uv)^{a_i+1}-1}.\]
\end{description}

\noindent \textbf{1.5. Theorem.} \textsl{Let $X$ be a} projective \textsl{log terminal variety of dimension $d$. Then $E_{st}(X;u,v)=(uv)^dE_{st}(X;u^{-1},v^{-1})$.}\\

This is an easy consequence of the alternative expression for $E_{st}(X)$ above, using the similar equality (2) for Hodge-Deligne polynomials of the smooth projective varieties $D_I$ \cite[Theorem 3.7]{Batyrev}.\\

\noindent \textbf{1.6.} From now on we consider only Gorenstein varieties $X$. Since then all discrepancies are integers, we have that the conditions canonical and log terminal are equivalent. 

\textsl{Assume} now that for a Gorenstein canonical projective $X$ of dimension $d$ the rational function $E_{st}(X;u,v)$ is in fact a polynomial $\sum_{p,q} b_{p,q}u^pv^q$. (An intuitive formulation of this assumption is perhaps that $X$ is close to admitting a crepant log resolution). Then Batyrev defined the \textsl{stringy Hodge numbers} of $X$ as $h^{p,q}_{st}(X):=(-1)^{p+q}b_{p,q}$. Note as a motivation that
\begin{description}
\item[(a)] when $X$ is smooth we have $h_{st}^{p,q}(X)=h^{p,q}(X)$ since then $E_{st}(X)=H(X)$,
\item[(b)] by Remark 1.4  (iii) we have $h_{st}^{p,q}(X)=h_{st}^{q,p}(X)$,
\item[(c)] Theorem 1.5 yields that $h^{p,q}_{st}(X)$ can be nonzero only for $0\leq p,q\leq d$ and that $h_{st}^{p,q}(X)=h_{st}^{d-p,d-q}(X)$,
\item[(d)] looking at the expression for $E_{st}(X)$ in Remark 1.4 (iv) we have clearly $h^{0,0}_{st}(X)=1$.
\end{description}

So the classical properties of Hodge numbers for smooth projective varieties are still satisfied by Batyrev's stringy Hodge numbers for Gorenstein canonical projective varieties with a polynomial stringy $E$-function. But there is one serious problem: it is not clear that they are nonnegative$\,$!\\

\noindent \textbf{Conjecture} \cite[Conjecture 3.10]{Batyrev}\textbf{.} \textsl{Let $X$ be a Gorenstein canonical projective variety. Assume that $E_{st}(X;u,v)$ is a polynomial. Then all stringy Hodge numbers $h^{p,q}_{st}(X)$ are nonnegative.}\\

\noindent \textbf{1.7.} The conjecture is trivially true for varieties that admit a crepant resolution (see Remark 1.4 (i)). This is the case for all canonical surfaces, which are exactly those with two-dimensional $A$-$D$-$E$ singularities. (In \cite{Veys2} there is related work on $\mathbb{Q}$-Gorenstein log terminal surfaces.)

In higher dimensions other results are rare. The first author settled the case of $A$-$D$-$E$ singularities of arbitrary dimension. In \cite{schepers} he was able to prove quite elegant concrete formulae for their stringy $E$-function, obtaining the following as a corollary.\\

\noindent \textbf{Theorem.} \textsl{Let $X$ be a projective variety of dimension at least 3 with at most $A$-$D$-$E$ singularities. Its stringy $E$-function is a polynomial if and only if $\dim X=3$ and all singularities are of type $A_n$ ($n$ odd) and/or $D_n$ ($n$ even). In that case, the stringy Hodge numbers of $X$ are nonnegative.}\\

\noindent \textbf{1.8.} In this paper we prove the following theorem (see Theorem 3.1).\\

\noindent \textbf{Theorem.} \textsl{Let $Y$ be a projective variety of dimension $d\geq 3$ with at most isolated Gorenstein canonical singularities. Assume that $Y$ admits a log resolution $f:X \to Y$ with all discrepancy coefficients of exceptional components $> \lfloor \frac{d-4}{2} \rfloor$. Write $E_{st}(Y;u,v)$ as a power series $\sum_{i,j\geq 0} b_{i,j} u^iv^j$. Then $(-1)^{i+j}b_{i,j} \geq 0$ for $i+j \leq d$.}\\

In view of 1.6 (c) this implies nonnegativity of the stringy Hodge numbers for such varieties $Y$ if they have a polynomial stringy $E$-function. As a corollary one can prove the analogous statement for arbitrary projective threefolds with Gorenstein canonical singularities (Corollary 3.2), and thus one obtains a proof of Batyrev's conjecture for threefolds in full generality. It is interesting that the somewhat mysterious polynomial condition on the stringy $E$-function in Batyrev's conjecture can be dropped in the cases of the theorem and the corollary. This suggests to use the numbers $b_{i,j}$ in these cases to generalize Batyrev's stringy Hodge numbers (Remark 3.3 (ii)). In addition we also give an explicit geometric description of these numbers $b_{i,j}$ (and thus in particular for the stringy Hodge numbers if they are defined) as dimensions or sums of dimensions of Hodge components of certain cohomology groups (Propositions 3.4 and 3.5).\\

\noindent \textbf{Acknowledgement:} We wish to thank Joost van Hamel for helpful discussions about the proof of Theorem 3.1.\\

\section{Preliminaries}

\noindent \textbf{2.1.} We will use the following result of de Cataldo and Migliorini in a crucial way. By $H^{\bullet}(\,\cdot\,)$ we still mean singular cohomology with coefficients in $\mathbb{C}$.\\

\noindent \textbf{Theorem.} \textsl{Let $Y$ be a projective variety of dimension $n$ with at most isolated singularities. Let $f: X \to Y$ be a resolution of singularities with $X$ projective, such that $f:f^{-1}(Y_{ns})\to Y_{ns}$ is an isomorphism, where $Y_{ns}$ is the nonsingular part of $Y$, and such that $f^{-1}(y)$ is a divisor for every singular point $y\in Y$. Denote by $D$ the total inverse image of the singular points. Then the map $H^i(X)\to H^i(D)$, induced by inclusion, is surjective for $i\geq n$.}\\

The proof of this result for $n=3$ and for one singular point is given in \cite[Theorem 2.3.4]{dCM2}. From this proof it is clear that it works for any dimension and any number of singular points, but for completeness we include the argument anyway. \\

\noindent \textbf{Proof:} We embed $Y$ in a projective space $\mathbb{P}^r$ and we take a generic hyperplane section $Y_s$ (so this hyperplane section is smooth and does not contain any singular point of $Y$). Consider the inverse image $X_s := f^{-1}(Y_s)$ and denote $Y\setminus Y_s$ by $Y_0$ and $X\setminus X_s$ by $X_0$. Let $\overline{Y_0}$ be the normalization of $Y_0$. Then $\overline{Y_0}$ still has isolated singularities and the map $f:X_0\to Y_0$ factors through $\overline{Y_0}$. Denote the map $X_0\to \overline{Y_0}$ by $\overline{f}$.
Because $\overline{Y_0}$ is affine, $H^i(\overline{Y_0}) =0 $ for $i> n$ (this was probably first proved by Kaup in \cite{Kaup}). The sheaves $R^k\overline{f}_*\mathbb{C}_{X_0}$ are skyscraper sheaves above the singular points for $k>0$, so they are flasque and $H^j(\overline{Y_0},R^k\overline{f}_*\mathbb{C}_{X_0})=0$ for $k,j> 0$. The Leray spectral sequence for $\overline{f}$ (given by $E_2^{p,q}=H^p(\overline{Y_0},R^q\overline{f}_*\mathbb{C}_{X_0})$ and abutting to $H^{\bullet}(X_0,\mathbb{C})$) gives then an isomorphism $H^i(X_0)\cong H^i(D)$ for $i> n$, because $E_2^{0,i}=E_{\infty}^{0,i} = H^i(D)$ and because $E_2^{p,q}=E_{\infty}^{p,q}=0$ for $p>0$ and $p+q=i$. There is also a surjection $H^n(X_0)\to H^n(D)$ because $E_2^{0,n}=E_{\infty}^{0,n} = H^n(D)$. Proposition (8.2.6) from \cite{Deligne2} states that the image of $H^i(X)$ in $H^i(D)$ is equal to the image of $H^i(X_0)$ in $H^i(D)$, and thus the result follows.\hfill $\blacksquare$\\

The following corollary is immediate from the theorem, see also \cite[Corollary 2.1.11]{dCM1}. It was already proved by Steenbrink in \cite[Corollary (1.12)]{Steenbrink} for the case where $D$ is a divisor with smooth components and normal crossings.\\

\noindent \textbf{Corollary.} \textsl{The a priori mixed Hodge structure on $H^i(D)$ is pure of weight $i$ for $i\geq n$.}\\

\noindent \textbf{2.2.} We also need the construction of the mixed Hodge structure on the cohomology of an algebraic set $D$ with smooth projective irreducible components and normal crossings. This can be found in \cite[p.149-156]{KulKur}. The construction is as follows. Denote the irreducible components of $D$ by $D_i,i\in I=\{1,\ldots, \alpha\}$ and put for $j\geq 0$
\[ D^{(j)}:= \coprod_{\substack{J\subset I\\ |J|=j+1}} D_J, \]
where $D_J=\cap_{i\in J} D_i$. So all $D^{(k)}$ are smooth and projective, and we write $A^i(D^{(k)})$ for the $C^{\infty}$-differential $i$-forms with values in $\mathbb{C}$ on $D^{(k)}$. The inclusion map $D^{(k)}\hookrightarrow D^{(k-1)}$ defined by mapping $D_{i_1}\cap \ldots\cap D_{i_{k+1}}$ into $D_{i_1}\cap \cdots\cap D_{i_{l-1}} \cap D_{i_{l+1}} \cap \ldots \cap D_{i_{k+1}}$ is denoted by $\delta_l^{(k)}$. 
There exists a spectral sequence $\{E_r,d_r\}$, abutting to $H^{\bullet}(D)$ and degenerating at the $E_2$ level, with $E_0^{p,q}=A^q(D^{(p)})$ and $d_0: E_0^{p,q}\to E_0^{p,q+1}$ given by differentiation of forms. Then $E_1^{p,q}$ is $H^q(D^{(p)})$; and $d_1: E_1^{p,q}\to E_1^{p+1,q}$ is defined by $\sum_{l=1}^{p+2} (-1)^l (\delta^{(p+1)}_l)^*$. There is a pure Hodge structure of weight $q$ on $E_1^{p,q}$, and $d_1$ is a morphism of Hodge structures, so $E_2^{p,q}$ inherits a pure Hodge structure and this provides $H^{\bullet}(D)$ with a mixed Hodge structure.\\

\section{Main results and corollaries} 

\noindent \textbf{3.1.} \textbf{Theorem.} \textsl{Let $Y$ be a projective variety of dimension $d\geq 3$ with at most isolated Gorenstein canonical singularities. Assume that $Y$ admits a log resolution $f:X \to Y$ with all discrepancy coefficients of exceptional components $> \lfloor \frac{d-4}{2} \rfloor$. Write $E_{st}(Y;u,v)$ as a power series $\sum_{i,j\geq 0} b_{i,j} u^iv^j$. Then $(-1)^{i+j}b_{i,j} \geq 0$ for $i+j \leq d$.}\\

\noindent \textbf{Remark.} 
\begin{description}
\item[(i)] In particular, if the stringy Hodge numbers of $Y$ are defined, then they are nonnegative in this case, in view of 1.6 (c). 
\item[(ii)] It is not very difficult to see that, if satisfied, then the condition on the discrepancy coefficients does not depend on the particular log resolution.
\item[(iii)] For $d=3$ there is no extra condition on the discrepancy coefficients (since we already assume that $Y$ is canonical) and for $d=4$ or 5 the condition means that $Y$ has (isolated Gorenstein) terminal singularities.
\end{description}

${}$

\noindent \textbf{Proof:} Let us first handle the case where $d\geq 4$. Choose a log resolution $f:X\to Y$ that is an isomorphism when restricted to $f^{-1}(Y_{ns})$ and denote the irreducible exceptional components by $D_l,l\in I=\{1,\ldots,\alpha\}$. Denote the discrepancy coefficient of $D_l$ by $a_l$ and put $D:= \cup_{l\in I} D_l$. Consider the alternative formula 1.4 (iv) for the stringy $E$-function:
\[E_{st}(Y) = \sum_{J\subset I} H(D_J; u, v) \prod_{l\in J} \frac{uv-(uv)^{a_l+1}}{(uv)^{a_l+1}-1},\]
and write this as a power series $\sum_{i,j\geq 0} b_{i,j}u^iv^j$. It suffices to prove that $(-1)^{i+j}b_{i,j}$ $\geq 0$ for $i\geq j$ and $i+j\leq d$. We denote the Hodge-Deligne polynomial of $X$ by $\sum_{i,j} c_{i,j}u^iv^j$ and for $J\subset I, J\neq \emptyset$ we denote the Hodge-Deligne polynomial of $D_J$ by $\sum_{i,j} c_{i,j}^{J}u^iv^j$.
The power series development of $\frac{uv-(uv)^{a_l+1}}{(uv)^{a_l+1}-1}$ for $a_l>0$ is equal to 
\begin{eqnarray*} 
& &(uv-(uv)^{a_l+1})(-1-(uv)^{a_l+1}-(uv)^{2a_l+2}-(uv)^{3a_l+3}- \cdots)\\ & & \quad = -uv + (uv)^{a_l+1} - (uv)^{a_l+2} + (uv)^{2a_l+2} -(uv)^{2a_l+3} + \cdots. \end{eqnarray*}
Since we assume the given lower bound on the $a_l$'s, we can write $b_{i,j}$ for $i+j\leq d$ and $i\geq j$ as
\[ b_{i,j}= c_{i,j} + \sum_{k=1}^{j} (-1)^k \sum_{\substack{J\subset I\\ |J|=k}} c_{i-k,j-k}^{J} + R_{i,j},      \]
where 
\[ R_{i,j} = \left\{\begin{array}{ll} \sum_{\substack{l=1,\ldots,\alpha\\ a_l=d/2-1}} c_{0,0}^{\{l\}} & \quad \text{if } d \text{ is even and } i=j=\frac{d}{2},\\ \sum_{\substack{l=1,\ldots,\alpha\\ a_l=(d-3)/2}} c_{0,0}^{\{l\}} &\quad \text{if } d \text{ is odd and } i=j=\frac{d-1}{2}, \\ \sum_{\substack{l=1,\ldots,\alpha\\ a_l=(d-3)/2}} c_{1,0}^{\{l\}} &\quad \text{if } d \text{ is odd and } i=\frac{d+1}{2}, j=\frac{d-1}{2},\\ 0 & \quad \text{otherwise}.\end{array}\right.\]
In any case, $R_{i,j}$ has got the right sign, so it suffices to study the rest of the formula for $b_{i,j}$. We rewrite it as follows, using the symmetry of usual Hodge numbers:
\[c_{d-i,d-j} + \sum_{k=1}^{j} (-1)^k \sum_{\substack{J\subset I\\ |J|=k}} c_{d-i,d-j}^{J}. \]
Set $q=2d-i-j$ and look at the $\{E_1^{p,q},d_1\}$-term of the spectral sequence from 2.2, for all $p$. This gives a complex:
\[ H^{q}(D^{(0)}) \to  H^{q}(D^{(1)}) \to \cdots \to  H^{q}(D^{(i-1)}) \to 0.\]
Here $H^q(D^{(i)})=0$ since $D^{(i)}$ has dimension $d-i-1$. The cohomology of this complex is given by $E_2^{0,q}$ to $E_2^{i-1,q}$. These are subquotients of $H^q(D)$ to $H^{q+i-1}(D)$ respectively. By purity of weight $r$ of $H^r(D)$ for $r > q$ (Corollary 2.1), $E_2^{1,q}$ up to $E_2^{i-1,q}$ must be zero$\,$! This means that the complex is exact. The purity of weight $q$ of $H^q(D)$ itself implies that $H^q(D) = \ker(H^q(D^{(0)})\to H^q(D^{(1)}))$. So we get an exact sequence where all arrows are morphisms of pure Hodge structures:
\[0 \to H^q(D) \to H^{q}(D^{(0)}) \to \cdots \to  H^{q}(D^{(i-1)}) \to 0.\]
We apply the exact functor $H^{d-i,d-j}$ and note that $H^{d-i,d-j}(H^q(D^{(s)}))= 0$ for $s=j,\ldots,i-1$ since $D^{(s)}$ has dimension $d-s-1$. Counting the dimensions of the resulting exact sequence of vector spaces gives
\[ \dim H^{d-i,d-j}(H^q(D)) = \sum_{k=1}^{j} (-1)^{k+1} \sum_{\substack{J\subset I\\ |J|=k}} (-1)^{2d-i-j}c_{d-i,d-j}^{J}. \] The surjectivity of the morphism of Hodge structures $H^q(X) \to H^q(D)$ from Theorem 2.1 translates by applying $H^{d-i,d-j}$ to 
\[ (-1)^{2d-i-j}c_{d-i,d-j} \geq \sum_{k=1}^{j} (-1)^{k+1} \sum_{\substack{J\subset I\\ |J|=k}} (-1)^{2d-i-j}c_{d-i,d-j}^{J},\]
and thus $(-1)^{2d-i-j}b_{i,j}=(-1)^{i+j}b_{i,j}\geq 0$. \\

For $d=3$ it can happen that some $a_l$ are zero. We can express the $b_{i,j}$ as follows:
\[ \left\{\begin{array}{l} b_{0,0}=c_{0,0}=1 \geq 0,\\ b_{1,0}=c_{1,0} \leq 0,\\ b_{2,0} = c_{2,0}  \geq 0,\\ b_{3,0}=c_{3,0}\leq 0, \\ b_{1,1} = c_{1,1}-\sum_{\substack{l=1,\ldots,\alpha\\ a_l\neq 0}} c_{0,0}^{\{l\}} = c_{2,2} - \sum_{\substack{l=1,\ldots,\alpha\\ a_k\neq 0}} c_{2,2}^{\{l\}} , \\ b_{2,1}=c_{2,1}-\sum_{\substack{l=1,\ldots,\alpha\\ a_l\neq 0}} c_{1,0}^{\{l\}}=c_{2,1}-\sum_{\substack{l=1,\ldots,\alpha\\ a_l\neq 0}} c_{2,1}^{\{l\}}.   \end{array}\right. \]
An analogous reasoning as above gives here that $H^4(D)=\oplus_{l=1}^{\alpha} H^4(D_l)$ and $H^3(D)=\oplus_{l=1}^{\alpha} H^3(D_l)$. Using the surjectivity of $H^i(X)\to H^i(D)$ for $i=3,4$ and applying the functors $H^{2,2}$, respectively $H^{2,1}$, immediately gives that $b_{1,1} \geq 0$ and $b_{2,1}\leq 0$. \hfill $\blacksquare$\\

\noindent \textbf{Example.} In particular, the above theorem and Batyrev's conjecture are true for fourfolds and fivefolds with isolated Gorenstein terminal singularities.\\ 

\noindent \textbf{3.2.} The problem in the following corollary is that the singularities can be non-isolated.\\

\noindent \textbf{Corollary.} \textsl{Let $Y$ be a projective threefold with at most Gorenstein canonical singularities. Write $E_{st}(Y;u,v)$ as a power series $\sum_{i,j\geq 0} b_{i,j}u^iv^j$. Then $(-1)^{i+j}b_{i,j}\geq 0$ for $i+j \leq 3$. In particular, Batyrev's conjecture is true for threefolds.}\\

\noindent \textbf{Proof:} The main theorem of \cite{Reid1} states that there exists a projective variety $Z$ with terminal singularities and a projective birational morphism $g:Z \to Y$ that is crepant. It follows then from \cite[Theorem 3.12]{Batyrev} that $E_{st}(Z)=E_{st}(Y)$. Since terminal singularities in dimension 3 are isolated (see for example \cite[Corollary 4-6-6]{Matsuki}), the corollary follows immediately from Theorem 3.1. \hfill $\blacksquare$\\ 

\noindent \textbf{3.3.} \textbf{Remark.} 
\begin{description}
\item[(i)] It is not true in general that $(-1)^{i+j}b_{i,j}\geq 0$ for $0\leq i,j\leq d$, as is shown by the following example (same notation as in the theorem above). Consider the variety $Z'=\{x_1^3x_5+x_2^4+x_3^2x_5^2+x_4^2x_5^2=0\}\subset \mathbb{P}^4$. The singular locus of $Z'$ consists of the line $\{x_1=x_2=x_5=0\}$ and an $E_6$ singularity at the point $P=\{(0:0:0:0:1)\}$. We can resolve the singular line by four blow-ups (first in the singular line itself, then in a surface isomorphic to $\mathbb{P}^2$ and afterwards consecutively in two curves isomorphic to $\mathbb{P}^1$). We are then left with a variety $Z$ which has a unique $E_6$ singularity at $P$. The stringy $E$-function of $Z$ can be written as $H(Z\setminus \{P\};u,v)\,+\,$(contribution of the singular point). For the second term we refer to \cite[Theorem 5.1]{schepers}, it is equal to $1+\frac{(uv)^2(2(uv)^6-2(uv)^5+(uv)^4-(uv)^2+2uv-2)}{(uv)^7-1}$. The first term can be computed to be $(uv)^3+7(uv)^2+7uv$. As a power series, $E_{st}(Z)=1+7uv+9(uv)^2-(uv)^3+(uv)^4-(uv)^6+\cdots$ and thus $b_{3,3}<0$.
\item[(ii)] It is tempting to generalize Batyrev's definition of stringy Hodge numbers for \textsl{any} variety $Y$ of dimension $d$ as in the theorem or the corollary as follows: for $i+j\leq d$ define analogously $h_{st}^{i,j}(Y):=(-1)^{i+j}b_{i,j}$ and maybe just put $h_{st}^{d-i,d-j}(Y):=h_{st}^{i,j}(Y)$.
\item[(iii)] In Theorem 3.1 we proved in fact that the numbers $(-1)^{i+j}b_{i,j}$ for $i+j\leq d$ and $d\geq 4$ are given by (with the same notations as in the theorem) \vspace{-0.6cm}

{\small \[(-1)^{i+j}b_{i,j} = \dim \ker (H^{d-i,d-j}(H^{2d-i-j}(X))\to H^{d-i,d-j}(H^{2d-i-j}(D))) + S_{i,j}, \]}where \[ S_{i,j} = \left\{\begin{array}{ll} \sum_{\substack{l=1,\ldots,\alpha\\ a_l=d/2-1}} \dim H^{d-1,d-1}(H^{2d-2}(D_l)) & \ \text{if } d \text{ is even and } i=j=\frac{d}{2},\\ \sum_{\substack{l=1,\ldots,\alpha\\ a_l=(d-3)/2}} \dim H^{d-1,d-1}(H^{2d-2}(D_l)) & \begin{array}{l} \text{if } d \text{ is odd and }\\ \quad i=j=\frac{d-1}{2}, \end{array} \\ \sum_{\substack{l=1,\ldots,\alpha\\ a_l=(d-3)/2}} \dim H^{d-2,d-1}(H^{2d-3}(D_l)) & \begin{array}{l}  \text{if } d \text{ is odd and }  \\  \quad  \{i,j\}=\{\frac{d-1}{2},\frac{d+1}{2}\},\\ \end{array}  \\ 0 &\ \text{otherwise}.\end{array}\right.\]   
\end{description}

In the following proposition we give a more intrinsic explicit description for these numbers. The 3-dimensional case is handled in 3.5.

${}$

\noindent \textbf{3.4.} \textbf{Proposition.} \textsl{Let $Y$ be as in Theorem 3.1 and of dimension $d\geq 4$. Take a log resolution $f: X\to Y$ that is an isomorphism when restricted to $f^{-1}(Y_{ns})$, with $D_l$ $(l=1,\ldots,\alpha)$ the irreducible components of the exceptional locus and with $a_l$ the discrepancy of $D_l$. Write $E_{st}(Y)$ as a power series $\sum_{i,j\geq 0} b_{i,j}u^iv^j$. Then for $i+j\leq d$ and we have \[(-1)^{i+j}b_{i,j} = \dim H^{d-i,d-j}(H^{2d-i-j}(Y)) + S_{i,j}, \]  
where $S_{i,j}$ is as in Remark 3.3 (iii).}\\

\noindent \textbf{Proof:} By Remark 3.3 (iii) above, it is sufficient to prove that 
\[\dim \ker (H^{d-i,d-j}(H^{2d-i-j}(X))\to H^{d-i,d-j}(H^{2d-i-j}(D))) \] \[ = \dim H^{d-i,d-j}(H^{2d-i-j}(Y)) \] 
for $i+j\leq d$. Let $y_1,\ldots,y_s$ be the (isolated) singular points of $Y$ and take disjoint contractible neighbourhoods $U_k$ for the complex topology around the points $y_k$. From the Leray spectral sequence for $f: V_k = f^{-1}(U_k) \to U_k$ it follows that $H^m(V_k)=H^m(f^{-1}(y_k))$ for all $m$. Now we consider the long exact sequences for $(X,\cup_k V_k)$ and $(Y,\cup_k U_k)$. We get the following diagram with exact rows (we use that $H^m(\cup_k V_k)=H^m(D)$ and $H^m(\cup_k U_k)=0$ for $m> 0$, where $D$ is the total exceptional locus): 

\begin{center}
\begin{picture}(115,16)
\put(0,12){$\to$}  \put(6,12){$H^{d-1}(D)$}  \put(23,12){$\to$} \put(29,12){$H^d(X,\cup_kV_k)$} \put(53,12){$\to$} \put(59,12){$H^d(X)$} \put(73,12){$\to$} \put(79,12){$\cdots$} \put(86,12){$\to$} \put(92,12){$H^{2d}(X)$} \put(107,12){$\to$} \put(113,12){$0$} \put(0,0){$\to$}  \put(11,0){0}  \put(23,0){$\to$} \put(29,0){$H^d(Y,\cup_kU_k)$} \put(53,0){$\to$} \put(59,0){$H^d(Y)$} \put(73,0){$\to$} \put(79,0){$\cdots$} \put(86,0){$\to$} \put(92,0){$H^{2d}(Y)$} \put(107,0){$\to$} \put(113,0){$0.$}
\put(11,6){$\uparrow$} \put(40,6){$\uparrow$} \put(63,6){$\uparrow$} \put(97,6){$\uparrow$} \put(113,6){$\uparrow$}
\end{picture}
\end{center} 

\noindent By excision we have 
\begin{eqnarray*}
H^m(Y,\cup_kU_k)&\cong& H^m(Y\setminus \{y_1,\ldots,y_s\},\cup_kU_k\setminus\{y_1,\ldots, y_s\})\\ &\cong &H^m(X\setminus D,\cup_kV_k\setminus D)\\ &\cong &H^m(X,\cup_kV_k)\end{eqnarray*} 
for all $m$ and then the Mayer-Vietoris construction (see \cite[Lemma 6.6]{Rotman}) gives us a long exact sequence (note that $H^{2d-1}(D)=H^{2d}(D)=0$)
\[ \to H^{d-1}(D) \to H^d(Y)\to H^d(X)\to H^d(D) \to H^{d+1}(Y)\to \cdots \] 
\[ \cdots \to H^{2d-1}(X)\to 0 \to H^{2d}(Y)\to H^{2d}(X)\to 0.\]
Thanks to Theorem 2.1 this sequence splits as follows:
\[ \to H^{d-1}(D) \to H^d(Y)\to H^d(X)\to H^d(D) \to 0, \tag{$d$}\]
\[ 0\to H^{d+1}(Y)\to H^{d+1}(X) \to H^{d+1}(D) \to 0, \tag{$d+1$}\]
\[  \vdots \]
\[ 0\to H^{2d-2}(Y)\to H^{2d-2}(X) \to H^{2d-2}(D) \to 0, \tag{$2d-2$}\]
\[ 0\to H^{2d-1}(Y)\to H^{2d-1}(X)\to 0, \tag{$2d-1$}\]
\[ 0\to H^{2d}(Y)\to H^{2d}(X)\to 0.\tag{$2d$}\]
Now we apply the exact functor $H^{d-i,d-j}$ to sequence $(2d-i-j)$ and we use Remark 3.3 (iii) to deduce the result. We must note that $H^{d-i,d-j}(H^{d-1}(D))=0$ for $i+j=d$. This is \cite[Corollary 2 p.154]{KulKur}, it also follows immediately from the discussion in 2.2. \hfill $\blacksquare$\\

\noindent \textbf{Remark.} We see from the short exact sequences $(d+1)$ to $(2d)$ that $H^k(Y)$ carries a pure Hodge structure for $k>d$. From the above proof it is clear that this works for any projective $Y$ with at most isolated singularities. Steenbrink proved this even for complete varieties with only isolated singularities \cite[Theorem (1.13)]{Steenbrink}.\\

\noindent \textbf{Example.} For a projective $Y$ of dimension $d=4$ or 5 and with at most isolated Gorenstein terminal singularities we can write for $i+j\leq d$ (notations as in the proposition)
\[ (-1)^{i+j}b_{i,j} = \dim H^{d-i,d-j}(H^{2d-i-j}(Y)) + \sum_{\substack{1\leq l \leq \alpha\\ a_l = 1}} H^{d-i+1,d-j+1}(H^{2d-i-j+2}(D_l)).     \] Indeed, for $d=4$ the second term can only be nonzero for $i=j=2$ and for $d=5$ this is the case for $i=j=2$ or $\{i,j\}=\{2,3\}$.\\

\noindent \textbf{3.5.} \textbf{Proposition.} \textsl{Let $Y$ be a projective threefold with Gorenstein canonical singularities. Write the stringy $E$-function of $Y$ as a power series $\sum_{i,j\geq 0} b_{i,j}u^iv^j$. Take a partial crepant resolution $g:Z\to Y$ as in the proof of Corollary 3.2. Then for $i+j\leq 3$ we have
\[ (-1)^{i+j}b_{i,j} = \dim H^{3-i,3-j}(H^{6-i-j}(Z)). \]}

\vspace{-0,5cm}

\noindent \textbf{Proof:} Recall from 3.2 that $Z$ has terminal singularities (which are automatically isolated) and that $E_{st}(Z)=E_{st}(Y)$. For $Z$ the analogue of Remark 3.3 (iii) (with $S_{i,j}$ always zero) is valid as well and thus we can proceed exactly as in the proof of the previous proposition. \hfill $\blacksquare$ \\

\footnotesize{

$\phantom{some place}$

\noindent Jan Schepers, Universiteit Leiden, Mathematisch Instituut, Niels Bohrweg 1, 2333 CA Leiden, The Netherlands\\
\textit{E-mail}: jschepers@math.leidenuniv.nl\\

\noindent Willem Veys, Katholieke Universiteit Leuven, Departement Wiskunde, Celestijnenlaan 200B, 3001 Leuven,
Belgium\\
\textit{E-mail}: wim.veys@wis.kuleuven.be

\end{document}